\documentstyle[12pt]{article}
\def\ba{\begin{eqnarray}}
\def\ea{\end{eqnarray}}

\def\g{\hbox{\bf g}}
\def\h{\hbox{\bf h}}

\def\lb{\label}
\def\be{\begin{equation}}
\def\ee{\end{equation}}

\begin{document}

\vspace{1cm}

\begin{center}
{\Large \bf On quantization of $r$-matrices
for Belavin-Drinfeld Triples}

\vskip 1cm
{A.P. Isaev$^*$ and O.V.Ogievetsky$^{\dagger}$}
\\[2em]
{\small $^*$ \it Bogoliubov Laboratory of Theoretical Physics, JINR,
141980  Dubna,} \\ {\small \it Moscow region, Russia}
\vskip .4cm
{\small $^{\dagger}$ \it Center of Theoretical Physics, Luminy,
13288 Marseille, France} \\
{\small and} \\
{\small \it P. N. Lebedev Physical Institute, Theoretical Department,
Leninsky pr. 53, 117924 Moscow, Russia} \\
\end{center}

\vskip 1cm
\centerline{{\bf Abstract}}
\vskip .3cm
\noindent
We suggest a formula for quantum universal $R$-matrices
corresponding to quasitriangular classical $r$-matrices
classified by Belavin and Drinfeld for all simple
Lie algebras. The $R$-matrices are obtained by twisting the
standard universal $R$-matrix.

\newpage

\noindent
{\bf 1.}
Classical quasitriangular $r$-matrices for semisimple Lie algebras
are classified by Belavin-Drinfeld triples \cite{BeDr}.
The Belavin-Drinfeld triple $(\Gamma_1, \, \Gamma_2, \, \tau)$
for a simple Lie algebra
$\g = \g^{+} \oplus \h \oplus \g^{-}$ consists of the following data:
$\Gamma_1$, $\Gamma_2$
are subsets of the set $\Gamma$
of simple roots of the algebra $\g$ and $\tau$
is a one-to-one mapping: $\Gamma_1 \rightarrow \Gamma_2$
such that $< \tau(\alpha), \, \tau(\beta) >= < \alpha, \, \beta >$
and $\tau^k(\alpha) \neq \alpha$
for any $\alpha,\beta \in \Gamma_1$ and any natural $k$.
The corresponding quantum $R$-matrices should have the form
\be
\lb{t0.1}
R_{12} = F_{21} \, {\cal{R}}_{12} \, F_{12}^{-1} \; ,
\ee
where ${\cal{R}}$ is
the standard universal Drinfeld-Jimbo $R$-matrix
for the Lie algebra $\g$. The twisting operator satisfies
the cocycle equation
\be
\lb{t0.2}
F_{12} \, (\Delta \otimes id) \, F =
F_{23} \, (id \otimes \Delta) \, F  \; .
\ee
Therefore the problem of quantization
is reduced to the problem of finding the twisting operator
$F_{12}$ for each Belavin-Drinfeld triple.
In the present paper we suggest a formula
for the twisting operator $F_{12}$. We present
the twisting  operator in a factorized form
\be
\lb{t0.2a}
F_{12} =  F^{(N)}_{12} \cdot F^{(N-1)}_{12} \cdots
F^{(2)}_{12} \cdot F^{(1)}_{12} \cdot K \; ,
\ee
where the factors $F^{(k)}$ are
special canonical elements
defined by the powers of the one-to-one map $\tau$;
the operator $K$ belongs to $q^{\h \otimes \h}$.
A different formula for the operator $F_{12}$ was given
in \cite{ESS}. We shall say several words about the differences at the end of
the present paper.

The plan of our paper is as follows.

Our approach heavily uses the modified Cartan-Weyl
basis for $U_q(\g)$. The definition of the modified
simple root generators is contained in Section 2. In Section 3 we give
an interpretation of Belavin-Drinfeld triples in terms
of the modified basis. In Section 4 a modified
Cartan-Weyl basis is introduced. The twisting operator
$F_{12}$ is constructed in Section 5. Finally, in Section 6
several examples are presented.

Everywhere below we assume the deformation parameter $q$ to be generic
(not a root of unity).

\vspace{0.5cm}

\noindent
{\bf 2.}  {\it Modified basis for quantum universal enveloping algebras.}

\vspace{0.5cm}

Consider a quantum universal enveloping algebra $U_q(\g)$ with
relations (see e.g. \cite{ChPr})
$$
%\begin{array}{c}
[h_i , \, h_j ] =0 \; , \;\;
[h_i , \, e_j ] = a_{ij} \, e_j \; , \;\;
[h_i , \, f_j ] = - a_{ij} \, f_j \; , \;\;
%\end{array}
$$
\be
\lb{t1.1}
[e_i , \, f_j ] = \delta_{ij} \,
\frac{K_i - K_i^{-1}}{q^{d_i} - q^{-d_i}} \; ,
\ee
and Serre relations
\be
\lb{t1.1e}
\sum^{1-a_{ij}}_{k=0} (-1)^k
\left[
\begin{array}{c}
1-a_{ij} \\ {k}
\end{array}
\right]_{q^{d_i}} \,
(e_i)^k \, e_j \, (e_i)^{1-a_{ij}-k} = 0 \; ,
\ee
\be
\lb{t1.1f}
\sum^{1-a_{ij}}_{k=0} (-1)^k
\left[
\begin{array}{c}
1-a_{ij} \\ {k}
\end{array}
\right]_{q^{d_i}} \,
(f_i)^k \, f_j \, (f_i)^{1-a_{ij}-k} = 0 \; ,
\ee
where
$$
\left[
\begin{array}{c}
n \\ {k}
\end{array}
\right]_{q} = \frac{[n]_q !}{[k]_q ! \, [n-k]_q !} \; ,
\;\;\; [k]_q = \frac{q^k -q^{-k}}{q-q^{-1}} \; ,
$$
$a_{ij}$ is the Cartan matrix for $\g$, $K_i = q^{d_i h_i}$ and  $d_i$
are smallest positive integers (from the set $1,2,3$) such
that $d_i a_{ij} = a^{(s)}_{ij}$ is symmetric matrix.
The algebra $U_q(\g)$ is a Hopf algebra with the comultiplication
\be
\lb{t1.2}
\begin{array}{c}
\Delta(h_i) = h_i \otimes 1 + 1 \otimes h_i \; ,  \\[3mm]
\Delta(e_i) = e_i \otimes K_i + 1 \otimes e_i \; , \;\;\;
\Delta(f_i) = f_i \otimes 1 + K_i^{-1} \otimes f_i \; .
\end{array}
\ee
The antipode and the counit are
$$
\begin{array}{c}
S(h_i) = -h_i \; , \;\;\;
S(e_i) = - e_i \, K_i^{-1} \; , \;\;\;
S(f_i) = - K_i \, f_i \; , \\[3mm]
\epsilon(h_i) = \epsilon(e_i) = \epsilon(f_i) = 0 \; .
\end{array}
$$

Any operator
$K \in q^{\h \otimes \h}$,
\be
\lb{t1.3}
K = q^{( \sum_{ij} b_{ij} \, h_i \otimes h_j )} \; ,
\ee
for arbitrary numerical matrix $b_{ij}$,
obviously satisfies the cocycle equation (\ref{t0.2}),
\be
\lb{t1.4}
K_{12} \, (\Delta \otimes id) \, K =
K_{23} \, (id \otimes \Delta) \, K  \; .
\ee
Therefore one can twist the comultiplication by $K$:
\be
\lb{twistk}
\widetilde{\Delta}(a) := K \, \Delta (a) \, K^{-1} \; .
\ee

We change the basis in the algebra $U_q(\g)$ by
introducing new generators
\be
\lb{t1.5}
E_i = X_i \, e_i \; , \;\;\; F_i =  f_i \, Y_i \; ,
\ee
where $X_i = \exp (\sum_j x_{ij} \, h_j)$,
$Y_i = \exp(\sum_j y_{ij} \, h_j)$
and $x_{ij}$, $y_{ij}$ are some numerical matrices.
We require that the comultiplication (\ref{twistk})
for the new generators (\ref{t1.5}) has the following form:
\be
\lb{t1.6}
\begin{array}{c}
\widetilde{\Delta}(E_i) = K \, \Delta(E_i) \, K^{-1} =
 E_i \otimes R^{+}_i + 1 \otimes E_i \; ,  \\ \\
\widetilde{\Delta}(F_i) = K \, \Delta(F_i) \, K^{-1} =
F_i \otimes 1 + R_i^{-} \otimes F_i \; .
\end{array}
\ee
Equations (\ref{t1.6}) relate
operators $X_i$, $Y_i$ and $K$.

A comparison of (\ref{t1.2}) and (\ref{t1.6}) gives
$$
X_i = q^{- \sum_{mn} h_m b_{mn} a_{ni}}
\equiv q^{-(h \, b \, a)_i} \; , \;\;\;
Y_i = q^{ \sum_{mn} h_m b_{nm} a_{ni}}
\equiv q^{(h \, \overline{b} \, a)_i} \; ,
$$
and
\be
\lb{t1.7}
R_i^{\pm} = X_i \, K_i^{\pm 1} \, Y_i =
K_i^{\pm 1} \, q^{-(h(b-\overline{b})a)_i} \; , \;\;\;
R_i^{+} = K_i^{2} \, R_i^{-} \; ,
\ee
where $\overline{b}_{mn} = b_{nm}$ is the transposed matrix.

The relations (\ref{t1.1})
and Serre relations (\ref{t1.1e}), (\ref{t1.1f})
for the quantum algebra $U_q(g)$ in
terms of the new generators (\ref{t1.5}) take the form
\be
\lb{t1.8a}
[E_i , \, F_j] = \delta_{ij} \frac{R^{+}_i - R^{-}_i}{q^{d_i} - q^{-d_i}}
\; ,
\ee
\be
\lb{t1.8}
\begin{array}{c}
R^{\pm}_i \, E_j = q^{\pm  a^{(s)}_{ij} + A_{ij}} \, E_j \, R^{\pm}_i \; ,
\;\;
R^{\pm}_i \, F_j = q^{\mp a^{(s)}_{ij} - A_{ij}} \, F_j \, R^{\pm}_i \; ,
\end{array}
\ee
\be
\lb{t1.8e}
\sum^{1-a_{ij}}_{k=0} (-1)^k
\left[
\begin{array}{c}
1-a_{ij} \\ {k}
\end{array}
\right]_{q^{d_i}} \, q^{-k A_{ij}}
\, (E_i)^k \, E_j \, (E_i)^{1-a_{ij}-k} = 0 \; ,
\ee
\be
\lb{t1.8f}
\sum^{1-a_{ij}}_{k=0} (-1)^k \left[
\begin{array}{c}
1-a_{ij} \\ {k}
\end{array}
\right]_{q^{d_i}} \,
q^{k A_{ij}} \,
(F_i)^k \, F_j \, (F_i)^{1-a_{ij}-k} = 0 \; ,
\ee
with a skewsymmetric matrix $A_{ij} = (\overline{a} (b - \overline{b})a)_{ij}$.

In the sequel we shall use $q$-commutators:
$[A,B]_\mu := A \, B - \mu \, B \, A$.
Relations (\ref{t1.8e}), (\ref{t1.8f}) can be conveniently rewritten
in terms of $q$-commutators. For example,
for $a_{ij} = 0$  the relations $[e_i , e_j] = 0 = [f_i , f_j]$
are rewritten as
$$
\, [E_i ,  E_j  ]_{q^{A_{ij}}} = 0 \; , \;\;\;
[ F_i ,  F_j  ]_{q^{- A_{ij}}} = 0 \; ,
$$
while for $a_{ij} = -1$ we have
\be
\lb{t1.8S}
\begin{array}{c}
[[E_{i} , E_{j}]_{\mu} , E_i]_{\nu} = 0 =
[E_{j} , [E_{i} , E_{j}]_{\mu} ]_{\nu} \; , \\ \\ \,
[[F_{i} , F_{j}]_{\nu} , F_i]_{\mu} = 0 =
[F_{j} , [F_{i} , F_{j}]_{\nu}]_{\mu}  \; .
\end{array}
\ee
where
$\mu = q^{d_i + A_{i \, j}}$, $\nu = q^{d_i - A_{i \, j}}$.

\medskip
\noindent
{\bf Remark.} The modified basis for multiparametric twistings of
$U_q(\g)$ has been considered by T.Hodges \cite{Hodg}.

\vspace{0.5cm}

\noindent
{\bf 3.}  {\it Modified basis and
Belavin - Drinfeld triples.}

\vspace{0.5cm}

All the data from the Belavin-Drinfeld triple can be
conveniently interpreted in terms of the modified basis
for a suitable matrix $b_{ij}$:

\medskip
\noindent
{\bf Proposition.} {\it
Let $\Gamma$ be the set of simple roots of $\g$,
$\Gamma_1$ and $\Gamma_2$ subsets of $\Gamma$ and
$\tau$ a one-to-one mapping:
$\Gamma_1 \rightarrow \Gamma_2$. Then the following equations
for the matrix $b_{ij}$
\be
\lb{t1.9}
R^{+}_{\alpha_i} = R^{-}_{\tau(\alpha_i)}
\;\; \forall \; \alpha_i \in \Gamma_1
\ee
where $R^{\pm}_{\alpha_i} \equiv R^{\pm}_i$,
admit a solution if and only if the triple
$(\Gamma_1 , \, \Gamma_2 , \, \tau)$ is the Belavin-Drinfeld
triple.
}

\medskip
\noindent
{\bf Proof.} Assume that a solution of equation (\ref{t1.9}) exists.
 We then need to prove that the mapping $\tau$ satisfies conditions:
\ba
\lb{1}
&&   1)  \; {\rm for} \; {\rm any} \; \alpha \in \Gamma_1 \; {\rm there \; is \;
a \; natural} \; k \;
{\rm for \; which} \; \tau^k(\alpha) \in \!\!\!\!\! / \, \Gamma_1 \, ,
\ \ \ \ \ \ \ \ \ \ \ \ \ \ \
\\[3mm]
\lb{2}
&&   2) \; \rm{for \; any} \; \alpha, \; \beta \in \Gamma_1 \; , \;\;
<\tau(\alpha),\tau(\beta)> = <\alpha , \beta > \; .
\ea

The condition (\ref{1}) means that $\tau$ has no cycles: $\tau^k(\alpha) ) \neq \alpha$
for all $\alpha \in \Gamma_1$ and $k >0$.

Indeed, assume that $\tau$ has a cycle, $\tau^k(\alpha) = \alpha$
for some $\alpha \in \Gamma_1$ and a natural $k$.
Take a minimal $k$ with this property.

Then $R^{+}_{\tau^k(\alpha)} = R^{+}_{\alpha}$ and
equations (\ref{t1.7}) imply
\be
\lb{t1.9re}
\begin{array}{l}
R^{+}_{\alpha} = R^{-}_{\tau(\alpha)} =
K^{-2}_{\tau(\alpha)} \, R^{+}_{\tau(\alpha)} = \dots =
\left( \prod_{i=1}^{k} \,
K^{-2}_{\tau^i(\alpha)} \right) \, R^{+}_{\tau^k(\alpha)}  \\[3mm]
= \left( \prod_{i=1}^{k} \,
K^{-2}_{\tau^i(\alpha)} \right) \, R^{+}_{\alpha}\; .
\end{array}
\ee
Therefore
$$
\left( \prod_{i=1}^{k} \,
K^{-2}_{\tau^i(\alpha)} \right) = 1 \; ,
$$
which contradicts to the independence of generators
in the Cartan subalgebra of $U_q(\g)$. Thus,
$\tau^k(\alpha)$ never equals $\alpha$ which proves
the condition (\ref{1}).

\vspace{2mm}
To prove the condition (\ref{2}), note that
eq. (\ref{t1.9}) is equivalent to the following condition
on the skewsymmetric matrix
$A_{mn} = (\overline{a} (b - \overline{b}) a)_{mn}$
\be
\lb{t1.10}
A_{im} + A_{m \tau(i)} + a^{(s)}_{im} +
a^{(s)}_{\tau(i) m} =0 \; ,
\ee
where the
subscript $m$ runs over all simple roots while $i$
numerates only roots from $\Gamma_1$. Eq. (\ref{t1.10})
is obtained by commuting both sides of eq. (\ref{t1.9})
with $e_m$ (or $f_m$). Here it is important that $q$ is not a root
of unity.

For indices $i,m$ corresponding to roots
$\alpha_i,\alpha_m \in \Gamma_1$,
eq. (\ref{t1.10}) can be rewritten in the following three
equivalent forms
\ba
\lb{t1.11}
&& A_{i\tau(m)} + A_{\tau(m) \tau(i)} +
a^{(s)}_{i\tau(m)} + a^{(s)}_{\tau(i) \tau(m)} =0 \; , \\[2mm]
\lb{t1.11.1}
&& A_{mi} + A_{i \tau(m)} + a^{(s)}_{mi} +
a^{(s)}_{\tau(m) i} =0 \; , \\[2mm]
\lb{t1.11.2}
&& A_{m\tau(i)} + A_{\tau(i) \tau(m)} +
a^{(s)}_{m\tau(i)} + a^{(s)}_{\tau(m) \tau(i)} =0 \; ,
\ea
 The combinations (\ref{t1.10}) $+$ (\ref{t1.11.1})
and (\ref{t1.11}) $+$ (\ref{t1.11.2}) of the equations are,
respectively,
\ba
\lb{t1.12}
&&2a^{(s)}_{im} = - a^{(s)}_{\tau(m) i} - a^{(s)}_{\tau(i) m}
- A_{m\tau(i)} - A_{i \tau(m)} \; , \\[2mm]
\lb{t1.12.1}
&&2a^{(s)}_{\tau(i)\tau(m)} = - a^{(s)}_{i\tau(m)} - a^{(s)}_{m\tau(i)}
- A_{m\tau(i)} - A_{i\tau(m)} \; .
\ea
Therefore,
\be
\lb{t1.13}
a^{(s)}_{im} = a^{(s)}_{\tau(i)\tau(m)} \; ,
\ee
which is equivalent to the second condition
$<\tau(\alpha_i),\tau(\alpha_m)> = <\alpha_i , \alpha_m >$
for the Belavin-Drinfeld triple.
$\bullet$

\medskip
\noindent
{\bf Remark 1.} The difference of eqs. (\ref{t1.11}) and (\ref{t1.11.1})
gives the following relation on the matrix $A_{ij}$
\be
\lb{t1.13a}
A_{im} - A_{\tau(i)\tau(m)} =
a^{(s)}_{mi} + a^{(s)}_{\tau(m) i} - a^{(s)}_{i \tau(m)}
- a^{(s)}_{\tau(i) \tau(m)} = 0 \; .
\ee
This shows that the map $\tau$ does not change
the modified basis.

\medskip
\noindent
{\bf Remark 2.}
Consider two sequences of sets
\be
\lb{tt1.1}
\begin{array}{l}
\Gamma_1 = \Gamma_1^{(0)} \supset \Gamma_1^{(1)}
\supset \Gamma_1^{(2)} \dots \supset \Gamma_1^{(N)}
\supset \Gamma_1^{(N+1)} = \emptyset \; , \\ \\
\Gamma_2 = \Gamma_2^{(0)} \supset \Gamma_2^{(1)}
\supset \Gamma_2^{(2)} \dots \supset \Gamma_2^{(N)} \; ,
\end{array}
\ee
defined by
$$
\Gamma^{(k+1)}_1 = \Gamma^{(k)}_1 \cap \Gamma^{(k)}_2 \; , \;\;\;
\Gamma^{(k)}_1 \stackrel{\tau}{\longrightarrow} \Gamma^{(k)}_2 \; .
$$
We assume that the set $ \Gamma_1^{(N)}$ is not empty.
The number $N$ is called the degree of the triple
$(\Gamma_1, \, \Gamma_2, \, \tau)$.

Introduce a set $\widetilde{\Gamma}^{(k)}_1 =
\tau^{-k-1}(\Gamma^{(k)}_2) \in \Gamma_1$.
Then the mapping $\tau^k$: $\widetilde{\Gamma}^{(k-1)}_1
\stackrel{\tau^k}{\longrightarrow}
\Gamma^{(k-1)}_2 \neq \, \emptyset$ also defines a Belavin-Drinfeld
triple
\be
\lb{bdk}
(\widetilde{\Gamma}^{(k-1)}_1, \, \Gamma^{(k-1)}_2, \tau^k) \; .
\ee

\vspace{0.5cm}
\noindent
{\bf 4.} {\it Modified Cartan-Weyl
basis and normal order of roots.}
\vspace{0.5cm}

Let $\Delta_{+}$ be the system of all positive
roots of $\g$ with respect to $\Gamma$.
A construction of Cartan-Weyl basis in terms
of the modified generators $E_i$ and $F_i$
is analogous to the usual procedure
for $U_q(\g)$ (see \cite{KhoTo}).

Recall the notion of a normal (convex) order in $\Delta_+$: the set $\Delta_+$
is ordered normally if any root $\gamma$ which is a sum of roots
$\alpha$ and $\beta$ is placed between $\alpha$ and $\beta$.

We write $\alpha <\beta$ if the root $\alpha$ is located to the left
of the root $\beta$.
For $\alpha <\beta$, the interval between roots $\alpha$ and $\beta$ is denoted by
$\{ \alpha,\beta\}$.

Given a normal order in $\Delta_+$, the modified Cartan-Weyl basis is constructed
by the following inductive procedure. The generators for the simple
roots are already defined. For a composite root $\gamma$, take a minimal
interval $\{ \alpha ,\beta\}$, $\alpha <\beta$, with $\gamma =\alpha
+\beta$ ("minimal" means that there is no subinterval
 $\{ \widetilde{\alpha} ,\widetilde{\beta}\} \subset \{ \alpha ,\beta\}$
for which $\gamma = \widetilde{\alpha} +\widetilde{\beta}$).
Assume that generators $E_\alpha$, $E_\beta$, $F_\alpha$ and $F_\beta$
were defined at previous steps. Then generators $E_\gamma$ and $F_\gamma$ are defined by
\ba
\lb{cog1}
&&E_{\gamma} = [E_{\alpha}, \, E_{\beta}]_{\mu} =
E_{\alpha} \, E_{\beta} - \mu E_{\beta} \, E_{\alpha}  \; ,\\[2mm]
\lb{t1.19}
&&F_{\gamma} =[F_{\alpha}, \, F_{\beta}]_{\nu} =
F_{\alpha} \, F_{\beta} - \nu F_{\beta} \, F_{\alpha}  \; .
\ea
where
$$
\mu = q^{-<\alpha,\beta> + <\alpha, \, A \, \beta>}
\; , \;\;\;
\nu = q^{-<\alpha,\beta> - <\alpha, \, A \, \beta>}\
$$
and $A$ is the operator with the matrix $A_{ij}$:
$$
<\alpha_i, \, A \, \alpha_j> = A_{ij} \; .
$$
If there are several possible minimal intervals $\{ \alpha,\beta \}$
for which $\gamma =\alpha +\beta$, the definitions
(\ref{cog1})-(\ref{t1.19}) give proportional results.

\medskip\noindent
{\it Note.} For the case $A_{ij} =0$
the definition (\ref{cog1})-(\ref{t1.19}) of composite roots
does not coincide with the definition in \cite{KhoTo} since we
use the comultiplication (\ref{t1.2}) which is different from
the comultiplication in \cite{KhoTo}.

\vspace{0.5cm}
\noindent
{\bf 5.} {\it Twisting operators $F_{12}$ for Belavin-Drinfeld
triples.}
\vspace{0.5cm}

For a given simple Lie algebra $\g$ fix a normal order in $\Delta_+$.

We need the expression for the inverse of the universal $R$- matrix for the
algebra $U_q(\g)$:
\be
\lb{runiv}
{\cal R}^{-1} =
\vec{\prod}_{\beta \in \Delta_{+} } \;
\exp_{q_{\beta}} \left( - \lambda \, a_\beta \,
(e_{\beta} \otimes f_{\beta}) \right) \cdot K^{(0)}  \; ,
\ee
where $q_{\alpha} = q^{<\alpha,\alpha>}$, $\lambda = q-q^{-1}$
and $K^{(0)} \in q^{\h \otimes \h}$. The product in eq. (\ref{runiv})
is the ordered product corresponding to the chosen normal order of
roots. For precise values of the constants $a_\beta$ see \cite{R},
\cite{KhoTo}.
The function $\exp_q$ is the standard $q$-exponent,
\be
\lb{t0.4}
\exp_q(u) = \prod_{n=0}^\infty (1+(q-1)uq^n)^{-1}=
\sum_{k=0}^{\infty} \, \frac{u^k}{k_q!}
\; , \;\;\; k_q = \frac{q^k -1}{q-1} \; .
\ee

Let $(\Gamma_1, \, \Gamma_2 , \, \tau)$ be a  Belavin-Drinfeld triple of degree
$N$. Define elements $F^{(k)}$ by
\be
\lb{t2.1}
F^{(k)}_{12} =
\vec{\prod}_{\beta \in \Delta^{(k)}_{+} } \;
\exp_{q_{\beta}} \left( - \lambda \, a_\beta \,
(E_{\beta} \otimes F_{\tau^k(\beta)}) \right) \; ,
\ee
where in the ordered product we keep terms corresponding
to only those roots $\beta$ for which $\tau^k(\beta)$ is
defined (that is, the element $e_\beta$ belongs to the subalgebra
with generators from the subset $\widetilde{\Gamma}_1^{(k)}$
defined in (\ref{bdk})). This is reflected in the notation
$\beta \in \Delta^{(k)}_{+}$.

The expression (\ref{t2.1}) can be given a form
\be
\lb{t2.1.1}
F^{(k)}_{12} = (1 \otimes T^k) \left( K \, {\cal R}^{-1} \, (K^{(0)})^{-1} \,
K^{-1} \right) \; ,
\ee
where the operator $T$ on the elements $F_\beta$ is defined by
$T(F_\beta) = F_{\tau(\beta)}$ wherever $\tau(\beta)$ is defined;
$T(F_{\beta}) =0$ otherwise. The operator $K$ corresponds to the solution
of eqs. (\ref{t1.9}) for the given Belavin-Drinfeld triple.

\vspace{0.5cm}
\noindent
{\bf Theorem}. {\it For the quantum algebra $U_q(\g)$ and the
Belavin-Drinfeld triple
$(\Gamma_1, \, \Gamma_2 , \, \tau)$ of degree $N$
the universal twisting
element $F_{12}$ is
\be
\lb{t2.2}
F_{12} =  F^{(N)}_{12} \cdot F^{(N-1)}_{12} \cdots
F^{(2)}_{12} \cdot F^{(1)}_{12} \cdot K  \equiv
\widetilde{F}_{12} \cdot K \;
\ee
with the factors $F^{(k)}_{12}$ defined in (\ref{t2.1}).
}

\medskip
We sketch the proof shortly. It is based on explicit
formulas for the coproduct of elements $F^{(k)}_{12}$:
\be
\lb{t2.4}
(\widetilde{\Delta} \otimes id) F^{(k)}_{12} =
F^{(k)}_{23} \, (K^{(k)}_{23})^{-1} \,
F^{(k)}_{13} \, K^{(k)}_{23} \; ,
\ee
\be
\lb{t2.5}
(id \otimes \widetilde{\Delta} ) F^{(k)}_{12} =
F^{(k)}_{12} \, (K^{(k)}_{12})^{-1} \,
F^{(k)}_{13} \, K^{(k)}_{12} \; ,
\ee
for some elements $K^{(k)} \in q^{\h \otimes \h}$.
The comultiplication $\widetilde{\Delta}$ is twisted as in (\ref{twistk}).

Next, one can verify the following identities
\be
\lb{t2.8m}
\widetilde{F}^{(k)}_{23} \,
\widetilde{F}^{(k+m)}_{13} \, \widetilde{F}^{(m)}_{12}  =
\widetilde{F}^{(m)}_{12} \,
\widetilde{F}^{(k+m)}_{13} \, \widetilde{F}^{(k)}_{23} \; ,
\ee
where $\widetilde{F}^{(k)} = F^{(k)} \cdot (K^{(k)})^{-1}$.

With the help of (\ref{t2.4}) - (\ref{t2.8m})
it is straightforward to check the cocycle condition (\ref{t0.2}).

\medskip
\noindent
{\bf Remark 1.} Another expression for the twisting element $F$
was suggested in \cite{ESS}. The expression in \cite{ESS} has a
factorised form as well. However, the factors $F^{(i)}$ are different;
one of the differences is that each factor
in \cite{ESS} contains terms from $q^{\h \otimes \h}$.
In our expression (\ref{t2.2}) all terms from $q^{\h \otimes \h}$
are collected; the price is the appearance of the modified
basis.

\medskip
\noindent
{\bf Remark 2.} The element $F$ in (\ref{t2.2}) satisfies
the following analogue of the linear ABRR equation \cite{ABRR}:
\be
(1\otimes T) (F_{12}{\cal R}^{-1}(K^{(0)})^{-1}K^{-1})
=F_{12}K^{-1}\ .
\ee

\vspace{0.5cm}
\noindent
{\bf 6.} {\it Examples.}

\vspace{0.5cm}
\noindent
{\bf i) $U_q(sl(3))$ case} (see \cite{Hodg}).\\
Here we have only one nontrivial Belavin-Drinfeld triple:

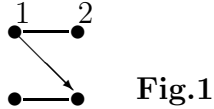
\unitlength=6mm
\begin{picture}(25,3)
\put(0.3,2){\circle*{0.3}}
\put(0.3,2.2){$1$}
\put(0.5,2){\line(1,0){1}}
\put(0.4,1.9){\vector(1,-1){1.2}}
\put(1.7,2){\circle*{0.3}}
\put(1.7,2.2){$2$}

\put(0.3,0.5){\circle*{0.3}}
\put(0.5,0.5){\line(1,0){1}}
\put(1.7,0.5){\circle*{0.3}}

\put(3,0.5){{\bf Fig.1}}
\end{picture}

This Cremmer-Gervais type triple has degree $1$ and the basic relations (\ref{t1.9})
which define this triple are reduced to one equation $R^{+}_1 = R^{-}_2$. The
antisymmetric matrix $A_{ij}$ is
\be
\lb{A}
A_{ij} = \delta_{i, j+1} - \delta_{j, i+1} \;\ ,
\ee
with $1\leq i,j\leq 2$.
The corresponding universal twisting
element (\ref{t2.2}) has the form
\be
\lb{t3.1}
F_{12} =   F^{(1)}_{12} \cdot K  = \exp_{q^2} (- \lambda \,
E_{1} \otimes F_{2})  \cdot K \; .
\ee

\vspace{0.5cm}
\noindent
{\bf ii) Cremmer-Gervais $U_q(sl(4))$ case.} \\
For this case the triple is given by the following diagram

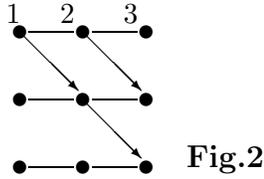
\unitlength=6mm
\begin{picture}(25,4)
\put(0.3,3){\circle*{0.3}}
\put(0,3.2){{\small $1$}}
\put(0.5,3){\line(1,0){1}}
\put(0.4,2.9){\vector(1,-1){1.2}}
\put(1.7,3){\circle*{0.3}}
\put(1.2,3.2){{\small $2$}}
\put(1.9,3){\line(1,0){1}}
\put(1.8,2.9){\vector(1,-1){1.2}}
\put(3.1,3){\circle*{0.3}}
\put(2.6,3.2){{\small $3$}}

\put(0.3,1.5){\circle*{0.3}}
\put(0.5,1.5){\line(1,0){1}}
\put(1.7,1.5){\circle*{0.3}}
\put(1.9,1.5){\line(1,0){1}}
\put(1.8,1.4){\vector(1,-1){1.2}}
\put(3.1,1.5){\circle*{0.3}}

\put(0.3,0){\circle*{0.3}}
\put(0.5,0){\line(1,0){1}}
\put(1.7,0){\circle*{0.3}}
\put(1.9,0){\line(1,0){1}}
\put(3.1,0){\circle*{0.3}}

\put(4,0){{\bf Fig.2}}

\end{picture}
\vspace{0.1cm}

It has degree 2. The basic relations (\ref{t1.9})
which define this triple are $R^{+}_1 = R^{-}_2$,  $R^{+}_2 = R^{-}_3$.
The matrix $A_{ij}$ is given by (\ref{A}), now with $1\leq i,j\leq 2$.
The corresponding universal twisting
element (\ref{t2.2}) has the form
\be
\lb{t3.3}
F_{12} =   F^{(2)}_{12} \cdot F^{(1)}_{12} \cdot K  \; ,
\ee
where
\be
\lb{t3.4a}
F^{(2)}_{12} = \exp_{q^2} (- \lambda \, E_{1} \otimes F_{3})  \; ,
\ee

\be
\lb{t3.4b}
F^{(1)}_{12} = \exp_{q^2} (-\lambda \, E_{1} \otimes F_{2})  \,
\exp_{q^2} ( q^{-1}\lambda \,
[E_{12}] \otimes [F_{23}]_{q^2} )  \,
\exp_{q^2} (- \lambda E_{2} \otimes F_{3})  \; .
\ee

\vskip 3mm
\noindent
Here $ [ E_{12} ] =E_1E_2 -E_2E_1$ and $ [ F_{23} ]_{q^2} =F_2F_3-q^2F_3F_2$.

\medskip
\noindent
{\bf Remark.}
One can directly check that (\ref{t3.1}), (\ref{t3.3}) obeys
the cocycle conditions (\ref{t0.2}). For (\ref{t3.1}) this check
requires only the basic equation for the $q$-exponent,
\be
\lb{t0.4a}
\exp_q (y) \, \exp_q(x) = \exp_q(x+y)  \;\;\;{\rm if} \;\;\; x \, y = q \, y \, x \; .
\ee
For (\ref{t3.3}) one needs two more quantum identities.
The first one is the famous pentagon
identity (see e.g. \cite{Kash} and references therein)
\be
\lb{t0.3}
\exp_q (u) \, \exp_q(v) = \exp_q (v) \, \exp_q([u,v]) \, \exp_q(u) \; ,
\ee
where the operators $u$ and $v$ satisfy the commutation
(Serre) relations
$$
u \, [u,v] = q \, [u,v] \, u \; , \;\;\;
v \, [u,v] = q^{-1} \, [u,v] \, v \; .
$$
The second identity is
\be
\lb{t0.5}
\exp_{q^2}(E) \, \exp_{q^2} (-R^{+}) \, \exp_{q^2}(F) =
\exp_{q^2}(F) \, \exp_{q^2} (-R^{-}) \, \exp_{q^2}(E) \; ,
\ee
where $E$, $F$ and $R^{\pm}$ generate the algebra
$$
[E,F] = (R^{+}-R^{-}) \; , \;\;\; [R^+ , \, R^-] = 0 \; ,
$$
$$
R^{\pm} \, E = q^{\pm 2} \, E \, R^{\pm} \; , \;\;\;
R^{\pm} \, F = q^{\mp 2} \, F \, R^{\pm} \; .
$$

\vspace{0.5cm}
\noindent
{\bf Acknowledgments:} We are indebted to
V.~Fock, S.~Khoroshkin, S.~Pakuliak, V.~Tolstoy
and especially P.N.~Pyatov
for valuable discussions and comments. 
One of the authors (API) thanks the Centre
de Physique Theorique in Luminy (Marseille University)
for the hospitality during his visit
in April 1999 where the considerable
part of the results of this paper has
been obtained. This work was supported
in part by the CNRS grant PICS No. 608 and the RFBR grant No.
98-01-2033.

%\section*{References}

\end{document}